# Logical Dependence of Physical Determinism on Set-theoretic Metatheory

Justin Clarke-Doane, Columbia University[1]

## 1. Introduction

Baroque questions of set-theoretic foundations are widely assumed to be irrelevant to physics. Let us call this the *insularity thesis*. In this article, I argue that this thesis is doubtful, especially once determinism is read the way it is actually used in physics. The dependence reaches into mainstream, foundational, and frontier physics, including statistical mechanics, discrete dynamics, canonicalization problems of the kind raised by gauge fixing, and the Kerr interior.

I begin with three analytic layer constructions. They illustrate how, even holding fixed a differential expression or definitional scheme, determinism verdicts for possible systems can differ between canonical extensions of ZFC — Gödel's Axiom of Constructibility (V=L) and ZFC plus large cardinal assumptions sufficient to imply Projective Determinacy (LC/PD). The differences track three aspects of determinism. One is *coherence*, through measurability assumptions built into weak formulations. Another is *uniqueness*, via hypotheses such as strict convexity whose truth can be tied to set-theoretic regularity predicates. The final one is *identity*, by means of a definition that does not pick out the same initial datum across metatheories. These illustrations are conceptual in the same way that familiar determinism thought experiments are.

The paper's main claim concerns a second layer of determinism that appears to have been largely overlooked in the literature. Determinism as actually used is meant to withstand admissible changes of gauge, mesh, coarse-graining, regulator, readout procedure, boundary prescription, and the like. It is also usually stated as a claim about generic or almost sure behavior. Making this precise forces one to ask whether the definable event sets are measurable or Baire and whether the associated selection problems admit regular solutions. Those *regularity* facts famously divide between V=L and LC/PD, even when pointwise membership facts are absolute.

I isolate a template under which robust determinism profiles land at the $\Sigma^1_2$ level, the first level of the projective hierarchy at which regularity (not membership) can fail in V=L while holding under LC/PD. So ZFC alone may already not decide whether such a profile is regular. One can go further, however. If a determinism profile is $\Sigma^1_2$-complete under continuous reductions, then it is universally measurable under LC/PD and not under V=L. The paper proves two unconditional theorems of this kind. First, a computable nearest-neighbor Ising Hamiltonian on $\mathbb{Z}^3$, with one fixed zero-temperature Glauber schedule, yields a $\Sigma^1_2$-complete tail readout profile in which only the initial microstate varies (Theorem 5.1). Second, at the Kerr Cauchy horizon, the problem of choosing canonical representatives for continuous extension germs in a fixed collar contains the eventual agreement relation $E_0$. So no universally measurable rule can solve it in ZFC, while V=L supplies a projectively definable rule and PD forbids any projective one (Theorem 6.2). A Shoenfield-style no go theorem separates these selection results from bare extension existence, which is analytic and admits no exact Borel calibration to constructibility (Theorem 6.1).

The Ising theorem is a coding theorem about one fixed Hamiltonian. It is not a theorem about a standard phase diagram under a physically distinguished disorder law. Similarly, the Kerr theorem canonicalizes codes of the same germ in a fixed presentation. Sections 5.5 and 6.5 clarifies these boundaries. It is open whether a distinguished measure, topology, or model class in an everyday physical theory reaches the hard part. I call the systematic investigation of such questions *reverse physics*, on rough analogy with Friedman's and Simpson's reverse mathematics. Given the entanglement of set theory and physics, any such theory would need to be relativized to foundational theories ("assuming ZFC + LC/PD, this robustness profile is regular, but it is not assuming ZFC + V=L"), or, as Quine (1951, 1990) argued, the search for new axioms to settle undecidables would be part and parcel to the search for new physical laws.

## 2. Set-theoretic background: V=L, LC/PD, and regularity

### 2.1 Two metatheories

ZFC is the standard foundational theory for mainstream mathematics. It leaves many questions about sets of reals undecided. Among the best-known candidates for settling such questions are Gödel's Axiom of Constructibility, V=L, and large-cardinal hypotheses strong enough to imply Projective Determinacy. I use "LC/PD" as shorthand for the latter package (Jech 2003; Kanamori 2009; Martin and Steel 1989; Maddy 1997). V=L asserts that every set is constructible. Philosophically, it is often read as a minimality principle: there are no sets beyond those generated by the constructible hierarchy. Most set theorists appear to reject it, though some defend it as mathematically (or, in Quine's case, physically) natural.[2] Large-cardinal hypotheses assert the existence of strong infinite cardinals. The existence of infinitely-many Woodin cardinals imply PD. Under PD, every projective set of reals is Lebesgue measurable, universally measurable, and has the Baire property (Kechris 1995; Moschovakis 2009; Martin and Steel 1989).

### 2.2 Projective complexity and coding

The projective hierarchy classifies sets of reals by definitional complexity. Starting with Borel sets, one obtains analytic sets, $\Sigma^1_1$, by continuous images or projections of Borel sets; coanalytic sets, $\Pi^1_1$, are their complements. Further alternations of real quantifiers produce $\Sigma^1_2$, $\Pi^1_2$, and the higher projective levels.

Physical state spaces are rarely literally $\mathbb{R}$ or $\omega^\omega$. The relevant spaces of fields, measures, histories, gauges, or initial data are often Polish or standard Borel spaces. To say that such an object is "coded by a real" is to fix a standard Borel presentation of the space. A "set of reals given by a code" is then the subset selected by the associated Borel or projective code. The coding is not innocuous under arbitrary changes of representation. Universal measurability is invariant under Borel isomorphism, not under every imaginable recoding. All claims below are made relative to fixed standard Borel presentations (Kechris 1995; Moschovakis 2009).

The primitive equations and constitutive data used in physics are often Borel. The rise in complexity typically occurs one level up, when the claim says that there exists a completion, selector, or continuation that works for every admissible implementation or diagnostic. That quantifier pattern, rather than the use of exotic coefficients by itself, is the main interface studied here.

### 2.3 Absoluteness and regularity

Shoenfield absoluteness insulates a large part of ordinary analysis from forcing and from passage to suitable inner models. In one standard form, $\Sigma^1_2$ and $\Pi^1_2$ sentences of second-order arithmetic, with parameters in the smaller model, are absolute between a universe and a transitive inner model containing all ordinals. In particular, membership in a fixed lightface $\Delta^1_2$ set is stable for reals common to the models (Shoenfield 1961; Jech 2003).

Shoenfield absoluteness does not support an indiscriminate comparison of every object described by the same words. V=L and a universe satisfying strong large cardinal hypotheses need not contain the same reals, and an intensional definition need not pick out the same extensional object in both. The important fact is that pointwise assertions can remain absolute across standard metatheoretic changes while *regularity properties of whole definable sets* do not. Whether a projective set is universally measurable or has the Baire property, and whether a projectively definable relation admits a selector of a specified regularity, can depend on the background axioms (Kechris 1995; Moschovakis 2009; Martin and Steel 1989).

### 2.4 A canonical $\Delta^1_2$ regularity failure in L

Fix a parameter-free lightface $\Delta^1_2$ formula which, in L, defines the canonical well-order $<_L$ of the reals of order type $\omega_1$. Fix also a Borel isomorphism $\iota:\omega^\omega\to(0,1)$ with a recursive Borel code. Let the same formula define $W\subseteq(0,1)^2$ by

$$W=\{(x,y):\iota^{-1}(x)<_L\iota^{-1}(y)\}.$$

In V=L, each horizontal section $W^y$ is countable, because every proper initial segment of the canonical order has countable order type, while each vertical section $W_x$ is co-countable. If W were Lebesgue measurable, Fubini's theorem would give measure 0 from the horizontal sections and measure 1 from the vertical sections. Hence W is not Lebesgue measurable and therefore not universally measurable.

The same sectional pattern rules out the Baire property in L, by the Kuratowski–Ulam theorem applied to both section families. The countable horizontal sections force W to be meagre, while the co-countable vertical sections force W to be comeagre. Appendix A.1 gives both applications. The construction can also be carried out directly on $2^\omega\times2^\omega$. Composing with a fixed recursive homeomorphism $2^\omega\times2^\omega\cong2^\omega$ gives a lightface $\Delta^1_2$ set $D\subseteq2^\omega$ which is not universally measurable and lacks the Baire property in L. This is the form used in Theorem 4.3.

Under LC/PD, the set defined by the same lightface $\Delta^1_2$ formula is projective and therefore universally measurable and has the Baire property. It need not still be a well-order. What is fixed across the comparison is the defining formula; what varies is the regularity, and possibly the extension, of the object it defines. This distinction will recur in the analytic layer examples.

## 3. Analytic determinism for possible systems: three aspects

Analytic layer determinism is the familiar one. It asks whether a law has solutions, whether they are unique, and whether the present state fixes a determinate evolution in the chosen function space. Standard treatments of PDEs and variational problems separate existence, uniqueness, and continuous dependence (Hadamard 1923; Earman 1986; Brezis 2011; Evans 2010).

The examples in this section are deliberately simple. They are not claims about actual materials, black holes, or finite laboratory systems. Like familiar determinism thought experiments, including Norton's Dome, their function is to test what follows from a stipulated mathematical law package rather than to describe a practicable experiment (Norton 2008). But unlike familiar determinism thought experiments, the present ones show how standard analytic notions can depend on the background set-theoretic metatheory which is almost universally ignored once coefficients, structural parameters, or initial data are fixed intensionally by projective definitions.

### 3.1 Coherence: weak formulations and measurability

Let $\Omega=(0,1)^2$ and consider the heat equation with a potential term,

$$\partial\square u=\Delta u-V(x)u,\quad u(0,\cdot)=u_0.$$

A standard weak formulation uses the bilinear form

$$a(u,v)=\int_\Omega(\nabla u\cdot\nabla v+Vuv)\,dx$$

on $H^1_0(\Omega)$. Measurability of V is part of the ordinary meaning of this expression. $V\in L^\infty(\Omega)$ is more than enough for the usual form methods.

Now let W be given by the fixed $\Delta^1_2$ formula of Section 2.4 and put $V=1_W$. Under LC/PD, W is measurable, so V is a bounded measurable potential and the weak formulation is coherent. Under V=L, W is not measurable, so the Lebesgue integral in the displayed bilinear form is not defined for all admissible u,v. Appendix A gives an explicit test function witnessing the failure.

One can alter the mathematics by using outer integrals, capacities, or a restricted test class. That produces a different formulation. The upshot is that the standard weak formulation attached to the fixed defining formula is not settled by ZFC.

### 3.2 Uniqueness: a regularity-sensitive structural parameter

Uniqueness theorems often depend on convexity, monotonicity, or Lipschitz hypotheses. These can be made to depend on a projective regularity statement without changing the elementary variational form.

Let $U\subseteq\mathbb{N}\times\mathbb{R}$ be a universal lightface $\Pi^1_2$ set, so every lightface $\Pi^1_2$ subset of $\mathbb{R}$ is a section $U_n$. Let

$$\psi\Leftrightarrow\ \forall\, n\in\mathbb{N}\ (U_n \text{ is Lebesgue measurable}),$$

and define $\gamma=-1$ if $\psi$ holds and $\gamma=+1$ otherwise. Consider

$$E_\gamma(z)=z^4-\gamma z^2.$$

Under LC/PD, every projective section is measurable, so $\gamma=-1$ and $E_\gamma$ is strictly convex with unique minimizer 0. In L, a lightface $\Delta^1_2$ nonmeasurable set occurs among the sections, so $\gamma=+1$ and the global minimizers are $\pm1/\sqrt{2}$. The same definitional rule yields a unique equilibrium in one metatheory and two in the other. The construction can be lifted to a coercive functional on a Hilbert space by applying the displayed energy to one mode and adding a positive quadratic term on its orthogonal complement (Dacorogna 2008).

### 3.3 Identity: when a definition does not fix the same datum

Finally, let $C\subseteq\mathbb{N}$ be

$$C=\{n:U_n \text{ is not Lebesgue measurable}\}.$$

Under LC/PD, $C=\varnothing$. In L, $C\neq\varnothing$. Choose smooth bumps $\varphi_n$ with pairwise disjoint supports and with

$$\|\varphi_n\|_{C^m}\leq 2^{-n}\ \text{for every } m\leq n,$$

and define

$$u_0(x)=\Sigma_{n\in\mathbb{N}}2^{-n}1_C(n)\varphi_n(x).$$

The series converges in $C^\infty$. Under LC/PD it is identically zero. But in L it is nonzero. So, a well-posed linear PDE has a unique solution for either datum, but the same definition does not identify the same initial state. The difference is one of identity rather than nonuniqueness.

### 3.4 Taking stock: relative possibility and infinite idealization

These examples establish a claim about mathematical law packages. They do not establish physical possibility *simpliciter*. A common proposal is explicitly relative. A history is physically possible relative to a set of laws when it is logically compatible with those laws (Maudlin 2007). To infer possibility under the actual laws, one must first identify the actual laws. The examples above are possible relative to the stipulated equations and definitions. But I have not argued that those are physically possible. There is a similar objection about finite systems. A finite lattice, a finite state space, or a finite family of admissible implementations is Borel-safe. The projective phenomena below require infinite spaces of data or choices. It is open to a critic to argue that a particular infinite idealization has been made incorrectly.

However, any such objection would cut across physical thought experiments wholesale. Nobody knows what the final laws will be, whether Norton's Dome or even GR spacetimes are physically possible (rather than merely relatively possible assuming the physical correctness of a merely effective theory). The criticism must thus identify which idealizations and quantifiers a physical theory licenses. It cannot be inferred merely from the fact that a foundational dependence appears, unless one builds the desired insularity into the criterion of acceptable idealization. (My own view is that no sense can be made of a distinguished notion of physical possibility, particularly because no sense can be made of a distinguished notion of logical possibility, and any notion of physical possibility would assume a notion of logical possibility (Clarke-Doane 2019, 2025). But nothing in what follows will depend on this point of view.)

The rest of the paper asks what happens when physically motivated robustness and selection demands do quantify over Polish spaces of completions and implementations – even if most do not carry high projective complexity.

## 4. The regularity layer: robustness, typicality, and canonicalization

### 4.1 Three regularity switches

At the regularity layer, the analogue of *coherence* is whether an event supports the probability or category language used to describe it. "Almost surely" presupposes measurability relative to a specified measure. A measure-independent assurance requires universal measurability. "Generically," in the Baire category sense, presupposes the Baire property.

The analogue of *uniqueness* concerns claims such as "for almost every input there is a unique admissible continuation." Even when existence and uniqueness are pointwise notions, the ensemble statement requires a regular event set.

The analogue of *identity* is canonicalization. Given a relation $R\subseteq X\times Y$ of inputs and admissible completions, a canonicalization demand asks for a selector $s:X\rightarrow Y$ with $(x,s(x))\in R$. A physical application will usually impose further requirements such as gauge compatibility, covariance, measurability, continuity on a large set, or stability under a specified class of representations. Mere choice is not enough (Kechris 1995; Moschovakis 2009).

*Definition 4.1.* A *regularity layer determinism profile* is a definable set or relation encoding when an input is deterministically well behaved only after a specified typicality, robustness, or canonicalization requirement has been imposed.

Note that the definition does not settle whether the resulting property deserves the word "determinism." A reader who reserves that word for pointwise existence and uniqueness may substitute "robustness profile" throughout. The mathematical issue is unchanged.

### 4.2 A universal-tail recipe

Let X be a Polish input space, Y a Polish space of completions or continuations, and $\mathbb{P}$ a Polish space of admissible policies, gauges, refinements, or diagnostics. Let

$$\mathrm{Good}(x,y,\rho,k,n)$$

be Borel and mean that at tolerance $2^{-k}$ and stage n, the pair (x,y) passes the diagnostic implemented by ρ. Define

$$\mathrm{UT}(x,y) \Leftrightarrow \forall\, \rho\in\mathbb{P}\ \forall\, k\in\mathbb{N}\ \exists\, N\ \forall\, n\geq N\ \mathrm{Good}(x,y,\rho,k,n),$$

and

$$\mathrm{DET}(x) \Leftrightarrow \exists\, y\in Y\ \mathrm{UT}(x,y).$$

*Proposition 4.2 (Projective upper bound).* If Good is Borel and $\mathbb{P}$ is Polish, then UT is $\Pi^1_1$ and DET is $\Sigma^1_2$.

*Proof.* For fixed (x,y,ρ), the quantifiers over k,N,n are countable and preserve Borelness. The remaining universal real quantifier over ρ makes UT coanalytic. Projecting over the real y makes DET $\Sigma^1_2$. ▮

This is an upper bound schema, obviously not a theorem that every robustness claim used by a physicist is $\Sigma^1_2$-complete. Compactness, uniqueness, convexity, or special PDE structure may collapse the complexity to Borel or analytic levels. A lower bound must be proved separately.

### 4.3 Completeness and forced regularity divergence

Let X be Polish. A set $A\subseteq X$ is $\Sigma^1_2$*-complete under continuous reductions* when A is $\Sigma^1_2$ and every $\Sigma^1_2$ subset of $2^\omega$ is the continuous preimage of A.

*Theorem 4.3 (Forced universal-measurability divergence).* Let X be Polish and let $A\subseteq X$ be $\Sigma^1_2$-complete under continuous reductions. Then:

(1) In ZFC + LC/PD, A is universally measurable and has the Baire property.

(2) In ZFC + V=L, A is not universally measurable.

*Proof.* Under PD, every projective set is universally measurable and has the Baire property. For the second clause, universal measurability is preserved under Borel preimages. If A were universally measurable in L, every $\Sigma^1_2$ set would be universally measurable there by continuous reducibility. Section 2.4 supplies a lightface $\Delta^1_2$ counterexample. ▮

There is no corresponding Baire conclusion in L. A $\Sigma^1_2$-complete set can be carried homeomorphically into a closed nowhere dense copy of Cantor space. The image remains complete but is meagre and hence has the Baire property.

To sharpen the contrast, it is useful to distinguish two kinds of metatheory sensitivity. A determinism profile set exhibits *can-diverge* sensitivity when its regularity properties are not provable in ZFC. This is already enough to show that determinism claims may require metatheoretic supplementation, and the standard physical practice of ignoring set theory is risky. A determinism profile set exhibits *must-diverge* sensitivity when it is regular in LC/PD and not universally measurable in V=L. This is the worst case scenario – and it is actual, as we will see.

### 4.4 Universal measurability and a distinguished physical measure

Universal measurability is stronger than measurability under one probability measure.

*Proposition 4.4.* If A is not universally measurable, then there is some Borel probability measure $\mu$ for which A is not measurable in the completion of $\mu$. Nothing follows about measurability under a separately specified measure $\nu$.

Thus Theorem 4.3 does not show that a physically distinguished chance measure fails to assign a probability to A. It shows that no measure independent guarantee is available in L. To turn the result into a claim about physical chance, one must identify the measure supplied by the physical theory and prove that the hard profile is nonmeasurable for that measure, or prove that no such distinguished measure is available. This is one of the principal tasks of reverse physics.

Universal measurability also should not be identified with full representation invariance. It is preserved under Borel isomorphism and secures measurability for every Borel probability on the fixed standard Borel space. It says nothing about arbitrary non-Borel changes of coding, symmetry invariance, robustness under perturbation, or the existence of a canonical selector.

### 4.5 A minimal canonicalization corollary

Let $A \subseteq X$ be projective and suppose a two-valued canonicalization rule is required to choose type 1 exactly on A. Such a rule is the indicator $1_A$. Under LC/PD it is universally measurable. If A is $\Sigma^1_2$-complete, then in L it is not. This elementary observation isolates the selection issue without conflating it with symmetry, covariance, or physical uniqueness. Whether a physical theory supplies the two candidate types and privileges this division is a separate question.

## 5. Case study I: statistical mechanics and a fixed-Hamiltonian coding theorem

### 5.1 Typicality, thermodynamic limits, and the selection layer

Statistical mechanics is a natural setting in which to look for the regularity layer. The central explanatory claims are rarely pointwise claims about one exact phase point. They concern typical microstates, thermodynamic limits, phase selection, stability of records, and the persistence of macroscopic verdicts under changes of volume sequence, boundary condition, coarse graining, or sampling procedure. In disordered systems, the difficulty is explicit. A sequence of finite volume Gibbs states need not converge, and different volume or boundary procedures can generate different limiting states. Metastates were introduced to represent this residual variation by probability measures on infinite volume Gibbs states (Newman and Stein 1997; Newman, Read, and Stein 2022).

Still, many ordinary statistical mechanical event sets do not lie high in the projective hierarchy. Cylinder events are Borel. Under standard continuity and compactness assumptions, many robust finite window events are analytic and so measurable in ZFC. However, a further selection problem can arise when the coarse specification does not determine a unique limiting object.

Let x code the interaction, temperature, disorder, and coarse preparation information. Let y code a proposed Gibbs state, metastate, or phase assignment. Let ρ range over admissible volume and boundary procedures. If

$$\mathrm{Good}_{SM}(x,y,\rho,k,n)$$

means that the k-local marginals at volume stage n under procedure ρ agree with y to tolerance $2^{-k}$, then the demand for a completion that works for every procedure has the form

$$\exists\ y\ \forall\ \rho\ \forall\ k\ \exists\ N\ \forall\ n\geq N\ \ \mathrm{Good}_{SM}(x,y,\rho,k,n).$$

Under a Borel presentation this is $\Sigma^1_2$. That is an upper bound. Whether any physically distinguished model or disorder ensemble realizes the bound is not settled by the quantifier count.

### 5.2 A conservative robust-trigger calculation

The same restraint is visible in a simple Mentaculus-style trigger event (Loewer 2020). Let $\Gamma_{PH}$ be a Past-Hypothesis macroregion, let $x\mapsto s_x$ be a continuous readout map over a compact time interval, and let $K_{\varepsilon,\eta}$ be a compact space of permitted window and threshold perturbations. Suppose the rulebook relation

$$R_{\varepsilon,\eta}(x,w,\kappa)$$

is closed. The set of pairs (x,w) for which the rulebook succeeds for every $\kappa\in K_{\varepsilon,\eta}$ is then closed, and the set $A(\varepsilon,\eta)$ of microstates admitting such a window is analytic. If $(\varepsilon_n,\eta_n)\downarrow(0,0)$, the fully robust event

$$G=\bigcap_n A(\varepsilon_n,\eta_n)$$

is still analytic because the witnesses $w_n$ can be packaged into one real. So, G is measurable and has the Baire property in ZFC. A single operationally robust event need not be metatheoretically sensitive. The projection from a coanalytic completion relation, rather than robustness language by itself, is what produces the $\Sigma^1_2$ jump.

### 5.3 The fixed system

Let us now give an exact lower bound. The simple construction is useful because it keeps the physical law fixed while locating the logical work.

Partition $\mathbb{Z}^3$ computably into nearest neighbor dimers. For example, pair $(2j,k,\ell)$ with $(2j+1,k,\ell)$ for every $(j,k,\ell)\in\mathbb{Z}^3$. Put coupling J=1 on each dimer edge and coupling 0 on every other nearest-neighbor edge. There is no external field. Formally,

$$H(\sigma)=-\Sigma_{\{u,v\}\in D}\sigma(u)\sigma(v),$$

where D is the fixed dimer set and $\sigma\in\{-1,+1\}^{\mathbb{Z}^3}$. The infinite sum need not converge; zero temperature single site dynamics uses only the local energy difference, which is well defined. At each site there is exactly one nonzero-coupled neighbor. The unique energy minimizing update sets the spin equal to that neighbor. Thus any configuration that is constant on every dimer is a stationary point of every asynchronous zero temperature update order.

Now fix one computable fair update schedule once and for all. It visits every lattice site infinitely often. Because the configurations used below are stationary, no freedom in the dynamics is needed. The second real quantifier will range over *readout paths*, not update laws.

Let $T_{\alpha,\beta}$, for $\alpha,\beta\in\omega^{\omega}$, be a standard universal family of trees on $\omega^{<\omega}$ such that

$$WF_{\exists} = \{\alpha: \exists\beta\ (T_{\alpha,\beta} \text{ is well founded})\}$$

is $\Sigma^1_2$-complete, and membership of a finite node in $T_{\alpha,\beta}$ depends continuously on $(\alpha,\beta)$. Enumerate the dimers by the nodes $s\in\omega^{<\omega}$. For each $(\alpha,\beta)$ define an exact initial microstate $\sigma_{\alpha,\beta}$ by putting both spins of the dimer $D_s$ equal to -1 when $s\in T_{\alpha,\beta}$ and equal to +1 otherwise. Put every unassigned dimer, if the enumeration leaves any, in the +1 state. The map $(\alpha,\beta)\mapsto\sigma_{\alpha,\beta}$ is continuous in the product topology, and the fixed Glauber evolution leaves it unchanged.

Let $\mathbb{P}\subseteq(\omega^{<\omega})^{\omega}$ be the closed space of strict chains

$$\rho(0)\subsetneq\rho(1)\subsetneq\rho(2)\subsetneq\ldots.$$

Define the stage-n verdict by reading the dimer indexed by $\rho(n)$:

$$V(\alpha,\beta,\rho,n) = -1 \text{ if } \rho(n) \in T_{\alpha,\beta}, \text{ and } V(\alpha,\beta,\rho,n) = +1 \text{ if } \rho(n) \notin T_{\alpha,\beta}.$$

The readout relation is clopen in the finite data on which it depends.

### 5.4 The theorem

*Theorem 5.1 (Fixed-Hamiltonian Ising tail-readout completeness).* There are a single computable nearest-neighbor Ising Hamiltonian on $\mathbb{Z}^3$, a single computable fair zero temperature asynchronous Glauber update schedule, a continuous family of stationary initial microstates $(\sigma_{\alpha,\beta})$, and a closed Polish space $\mathbb{P}$ of coherent readout paths such that

$$DET_{reg}= \{\alpha: \exists\beta\ \forall\rho\in\mathbb{P}\ \exists N\ \forall n\geq N\ [V(\alpha,\beta,\rho,n)=+1]\}$$

is $\Sigma^1_2$-complete under continuous reductions.

*Proof.* Fix $(\alpha,\beta)$. If $T_{\alpha,\beta}$ is ill founded, choose an infinite branch b. The path $\rho_b(n)=b\restriction(n+1)$ belongs to $\mathbb{P}$ and remains in the tree at every stage, so the verdict is always -1. The universal-tail condition fails.

If $T_{\alpha,\beta}$ is well founded, let $\rho\in\mathbb{P}$. The chain cannot remain in the tree forever, for then its union would be an infinite branch. Hence $\rho(N)\notin T_{\alpha,\beta}$ for some N. Since a tree is downward closed, no extension of a node outside the tree belongs to the tree. Thus $V(\alpha,\beta,\rho,n)=+1$ for every $n\geq N$. Therefore

$$\forall\rho\in\mathbb{P}\ \exists N\ \forall n\geq N\ [V(\alpha,\beta,\rho,n)=+1]$$

holds exactly when $T_{\alpha,\beta}$ is well founded. Projecting over β gives $DET_{reg}=WF_{\exists}$. ∎

*Corollary 5.2.* Under LC/PD, $DET_{reg}$ is universally measurable and has the Baire property. Under V=L, it is not universally measurable.

The existential variable β does not change the external field or Hamiltonian. It indexes a microstate compatible with the coded input α. The dynamics, couplings, and readout rule are fixed.[3]

### 5.5 What the theorem does not prove

The theorem is a coding theorem. The dimerized lattice is a countable register bank. The evolution is stationary. The admissible readout paths represent branches through a tree. The theorem is not a result about spontaneous symmetry breaking, an Edwards-Anderson or Sherrington-Kirkpatrick disorder law, chaotic size dependence, or the selection of a thermodynamic phase by ordinary laboratory procedures. So, calling its verdict "phase selection" may obscure where the complexity comes from. Nor does the theorem answer the finite idealization objection (mentioned in Section 3.4).. Every *finite* truncation of the register has a finite state space and a tame verdict. The complexity appears in the infinite family and in quantification over all coherent paths. A critic may conclude that the corresponding idealization is physically unmotivated. But more would need to be said. Familiar local interaction syntax does not force low descriptive complexity. To obtain an insularity theorem one needs substantive restrictions on initial state classes, measures, schedules, or admissible macroscopic questions.

Finally, it should be flagged that $\beta$ is a completion only in a formal sense. The theorem gives a continuous family of exact microstates above an infinite code $\alpha$. It does not prove that one finitely specified preparation has all the $\sigma_{\alpha,\beta}$ as physically possible microscopic completions. Establishing that relation for a physically distinguished coarse graining is a missing bridge.

## 6. Case study II: Kerr, bare extendibility, and germ representatives

### 6.1 The current Kerr regularity picture

In general relativity, global hyperbolicity is the standard mathematical surrogate for deterministic evolution. Appropriate Cauchy data determine a maximal globally hyperbolic development, but the development may admit extensions beyond a Cauchy horizon. Strong Cosmic Censorship is a family of conjectures asserting generic inextendibility at a specified regularity threshold (Wald 1984; Dafermos 2014).

The subextremal Kerr interior is the central rotating case. Dafermos and Luk proved that the Lorentzian metric extends continuously across a nontrivial part of the Cauchy horizon for their class of dynamical vacuum data (Dafermos and Luk 2025). Linearized analyses indicate curvature blow up there (Sbierski 2023; Gurriaran 2025). Sbierski's 2026 note strengthens the linear instability result and supplies an input to the later nonlinear proof (Sbierski 2026b). Sbierski also proved that a suitable curvature blow up assumption yields local Lipschitz inextendibility without symmetry assumptions (Sbierski 2026a). Gurriaran proved curvature blow up and Lipschitz inextendibility near timelike infinity under a nonlinear Price-law-type assumption (Gurriaran 2026). Luk and Sbierski proved weak null singularity formation for a characteristic initial value problem with smooth data settling to a strictly rotating subextremal Kerr event horizon and satisfying appropriate upper and lower bounds. The resulting metric is continuously, but not Lipschitz, extendible (Luk and Sbierski 2026).

These results identify an important regularity window and prove consequential conditional or model class results within it. The descriptive set-theoretic results below are orthogonal to them.

### 6.2 Bare extension existence is analytically shallow

Fix a standard Borel space $\mathscr{D}$ of coded initial or characteristic data and a regularity class $\mathscr{R}$. Suppose candidate local extensions across a prescribed compact horizon crossing region are coded by a Polish space $\mathscr{Y}$, and that

$$R_{\mathscr{R}}(d,y)$$

is a Borel relation saying that y codes an $\mathscr{R}$-extension compatible with the development of d. Then

$$\mathrm{Ext}^{\mathscr{R}}(d) \Leftrightarrow \exists\, y\in\mathscr{Y}\, R_{\mathscr{R}}(d,y)$$

is analytic. This classification is conditional on the stated coding. Verifying that a particular weak Einstein solution concept has a Borel certificate relation may be analytically nontrivial. But once the presentation is available, the projective bound follows immediately.

The bound rules out one tempting route to a perfect L-axis switch.

*Theorem 6.1 (Shoenfield no-go for exact bare calibration).* Let $a\in\omega^\omega$ code a Borel map $x\mapsto d_x$ and an analytic extension predicate. Suppose

$$\mathrm{Ext}^{\mathscr{R}}(d_x) \Leftrightarrow x\in L[a]$$

for every $x\in\omega^\omega$. Then every real belongs to L[a]. Consequently, in any universe containing a real outside L[a], no such exact calibration exists.

*Proof.* Write the pullback predicate as $\exists\, z\, \varphi(x,z,a)$ with $\varphi$ arithmetical. If $x\in L[a]$ and the displayed equivalence holds in the universe, then $\exists\, z\, \varphi(x,z,a)$ holds there. This is a $\Sigma^1_1(x,a)$ statement, hence $\Sigma^1_2(x,a)$, so Shoenfield absoluteness gives a witness in L[a]. Therefore

$$L[a]\vDash\forall\, x\, \exists\, z\, \varphi(x,z,a).$$

The displayed sentence is $\Pi^1_2(a)$, so Shoenfield absoluteness transfers it to the universe. By the assumed equivalence, every real then lies in L[a]. ▮

The conclusion is about the reals, not the whole universe. It is $\mathbb{R}\subseteq L[a]$, not $V=L[a]$. It also does not say that bare extendibility is absolute between arbitrary models that disagree about the parameters or the coding space. It says that a Borel pullback of an analytic predicate cannot exactly be the full set of a-constructible reals in an anti-L[a] universe.

Theorem 6.1 rules out a literal equivalence between bare extension existence and constructibility at this projective level. Any positive comparison must concern a richer robustness or selection predicate, use a non-Borel data map, or abandon exact L[a] calibration.

### 6.3 A fixed-collar space of continuous metric perturbations

There is an unconditional selection theorem for a fixed collar of continuous metric codes. Fix a compact collar

$$N=K\times[0,1]$$

on the extension side of a Cauchy-horizon patch, with K compact and s=0 the horizon. Let $\bar{g}_0$ be one continuous Lorentzian extension metric on N agreeing with the given interior metric at s=0. Exact Kerr has such collars, and the Dafermos-Luk theorem supplies them in its dynamical setting. Choose a background Riemannian norm on symmetric two-tensors and let

$$B=C_0(N;\mathrm{Sym}^2T^*N)$$

be the separable Banach space of continuous perturbations h satisfying $h|_{K\times\{0\}}=0$, with the uniform norm. Compactness and nondegeneracy give $\varepsilon>0$ such that $\bar{g}_0+h$ is Lorentzian whenever $\|h\|_\infty\leq\varepsilon$. Let X be that closed ball. It is Polish. Define equality as germs in the fixed collar by

$$h\ E_{germ}\ h' \ \Leftrightarrow\ \exists\ m\geq 1\ \forall\ (p,s)\in K\times[0,1/m]\ \ h(p,s)=h'(p,s).$$

The relation is $F_\sigma$, hence Borel. A selector for $E_{germ}$ is a function $S:X\to X$ such that $S(h)\ E_{germ}\ h$ and $S(h)=S(h')$ whenever $h\ E_{germ}\ h'$. It chooses one coded representative of every equality-germ class.

### 6.4 The $E_0$ theorem

Let $E_0$ be eventual equality on $2^\omega$:

$$\gamma\ E_0\ \gamma' \ \Leftrightarrow\ \exists\ N\ \forall\ n\geq N\ \gamma(n)=\gamma'(n).$$

*Theorem 6.2 (Fixed-collar germ-representative dichotomy).* In the setting of Section 6.3:

(1) There is a continuous injection $\Theta:2^\omega\to X$ with compact image F such that

$$\gamma\ E_0\ \gamma' \ \Leftrightarrow\ \Theta(\gamma)\ E_{germ}\ \Theta(\gamma').$$

(2) In ZFC there is no universally measurable selector for $E_{germ}$ on X or on F, and there is no Baire measurable selector on X or on F.

(3) In ZFC + V=L there is a selector on X whose graph is projective.

(4) In ZFC + PD there is no selector on X, or on F, whose graph is projective.

*Proof sketch.* Choose pairwise disjoint marker balls $B_n$ whose upper s-coordinates tend to 0, and continuous tensor bumps $b_n$ supported in $B_n$, with $\|b_n\|_\infty$ decreasing rapidly to 0. Put

$$\Theta(\gamma)=\Sigma_n\gamma(n)b_n.$$

The sum is continuous, stays in X, and two perturbations agree on some sufficiently small collar exactly when their bit sequences eventually agree. This proves (1).

For (2), suppose S were universally measurable and put $r=S\circ\Theta$. The map r is $E_0$-invariant. Under the Bernoulli product measure, every measurable set invariant under finite coordinate flips has measure 0 or 1. It follows, by applying the zero-one law to preimages of a countable base of X, that r is almost everywhere constant. But the inverse image of one $E_{germ}$-class is contained in one countable $E_0$-class, which cannot be conull. The category argument is parallel. An $E_0$-invariant set with the Baire property is meagre or comeagre, so a Baire measurable invariant map is constant on a comeagre set, again impossible.

For (3), in V=L use the canonical constructible well-order of the reals and choose the least member of each Borel equivalence class. The graph is projective (no sharper pointclass claim is needed). For (4), a function with projective graph has projective preimages of open sets. Under PD those preimages are universally measurable and have the Baire property, contradicting (2). Full details are in Appendix F. ∎

Note that the theorem does not require the Glimm-Effros dichotomy, though it places the result in its natural setting. The embedded relation is exactly $E_0$, the canonical minimum nonsmooth countable Borel equivalence relation (Harrington, Kechris, and Louveau 1990; Kechris 1995).

### 6.5 Interpretation: what is selected

Theorem 6.2 is unconditional, but its physical scope is narrow. $E_{germ}$ identifies perturbations that are already the same continuous metric germ. The selector chooses one code for each germ class. It does not choose among inequivalent continuations of the interior spacetime. The theorem concerns representational canonicalization. It does not recover a unique future beyond the Cauchy horizon.

The fixed-collar qualification is also key. A diffeomorphism invariant version would compare extensions up to local isometry over the common past. The marker construction does not automatically embed $E_0$ into that larger equivalence relation. A diffeomorphism might permute or erase the chart readable markers. One would need geometric rigidity or invariant marker data. That is open.

Finally, the perturbations are continuous Lorentzian metrics satisfying the required boundary values. No vacuum Einstein equation is imposed on the extension side. This is appropriate for displaying the descriptive obstruction inherent in the bare $C^0$ metric-germ space, but it recommends against calling the theorem a SCC theorem. A vacuum constrained, weak solution, or isometry invariant counterpart would be substantially stronger.

Nevertheless, once a formalism represents weak continuations by a nonsmooth quotient, the instruction to choose a canonical representative adds mathematical content. V=L supplies a definable but nonregular choice. PD makes every projective rule regular and so rules out such a selector. The theorem does not decide whether this fixed presentation quotient is the physically right one.

## 7. Nearby foundational dependence without bespoke physical coding

The preceding theorems separate three claims. A projective upper bound can arise from the logical form of a robustness demand. A completeness result can be obtained by an explicit register construction. And an already given quotient problem can obstruct regular selection. The the general phenomenon is not an invention of the present coding.

### 7.1 Pitowsky spin models

Pitowsky proposed nonmeasurable hidden variable models for spin. Farah and Magidor proved that the existence of Pitowsky spin models is independent of ZFC (Farah and Magidor 2012). This is the closest established precedent that I am aware of in foundations of physics. It is also cautionary. The theorem concerns the mathematical availability of a particular hidden variable construction. It does not show that ordinary quantum experiments have different frequencies in different universes of set theory, and it does not overturn Bell's theorem under its standard probability assumptions.

The upshot is that "there is no such model" can contain more set theory than the quantum mechanical formalism advertises. Whether the model has an acceptable physical interpretation is a further question. Farah and Magidor explicitly distinguish the independence result from a direct physical interpretation. The same distinction governs the present paper.

### 7.2 Measurable equilibrium selection

Levy studies Borel parametrized families of games in which equilibria exist pointwise and asks whether they can be selected measurably. The existence of Lebesgue- or universally measurable equilibrium selectors is independent of ZFC, even under restrictions to zero sum games and with mixed strategies allowed (Levy 2025). This is the distinction between pointwise existence and regular global selection emphasized here. No universal tree is inserted into the statement of the economic problem. The foundational sensitivity belongs to the selection demand itself.

The importance of this result for present purposes is methodological. It blocks the assumption that whenever each parameter value has an admissible object, asking for a measurable family of such objects is merely routine analysis. The global selector can carry independent content.

### 7.3 Disintegration of measures

Backs, Jackson, and Mauldin connect Maharam's question about uniformly σ-finite disintegrations to Borel uniformization and show that the answer may depend on additional set theoretic axioms.  In particular, CH and V=L yield negative results and special $\Pi^1_1$ constructions (Backs, Jackson, and Mauldin 2015). Disintegration underlies conditional probability and ergodic theory. The lesson is that a natural global regularity question attached to ordinary measure theory may not be settled by ZFC.

### 7.4 Robust optimization under PD

Carassus and Ferhoune formulate robust utility maximization with projective data and use PD to obtain optimal strategies under their hypotheses. In separate work they develop closure, integration, and measurable selection results for projective functions, including ε-optimal selectors under PD (Carassus and Ferhoune 2024, 2025). These results show how the projective hierarchy can enter a mathematically motivated optimization problem without first being named by the modeler. They do not show that a market or physical system realizes every projective input, but what stronger regularity axioms do once such inputs and selection demands are admitted.

### 7.5 The comparative moral

These examples are more natural than the Ising register and less about physical determinism. There is foundational dependence in global selection and uniformization problems studied for independent reasons. There are also explicit coding theorems showing that physical looking local formalisms can realize maximal projective complexity. What remains scarce is a theorem that joins both virtues: a physically distinguished model and measure, together with an unengineered determinism profile whose regularity provably divides between serious extensions of ZFC.

The mathematics establishes a possibility and a research program. It does not yet establish that a familiar physical argument to date implicitly depends on a choice between V=L and LC/PD.

## 8. Objections and replies

I have argued that determinism is entangled with a canonical debate over axiom candidates extending standard set theory. First, at the analytic layer, the coherence, uniqueness, and identity of possible physical systems differ with V=L and LC/PD. Second, once determinism as it is actually invoked in physics is seen to involve a regularity layer, these three features can differ in physical systems too. How might one respond to these arguments? I can think of the following ways.

*Objection 1. Physicists rarely use projective level 2 or higher sets.*

This objection confuses explicit terminology with implicit structure. Consider the types of conditions physicists routinely impose: "There exists an admissible object y such that for every admissible choice ρ, y passes the robustness check." This has the logical form $(\exists y)(\forall \rho)\text{Good}(y,\rho)$, which is $\Sigma^1_2$ when y and ρ range over reals and Good is arithmetical or Borel. The mathematical structures that appear in mature physical theories naturally involve alternating quantifiers over infinite objects. The Kerr interior predicates of $C^0$-extendibility, finite flux, blow-up, choice-stability, and gauge independent certification are all core concepts from the GR literature (Wald 1984; Dafermos 2014). Similar remarks apply to metastate existence and robust trigger events (Newman and Stein 1997), and to gauge fixing independence (Gribov 1978; Singer 1978).

But let us suppose that this objection were correct as far as it went. What of philosophical interest would it show? Consider the mathematical case. Just as the MRDP Incompleteness reveal that seemingly elementary diophantine questions are independent of ZFC, this paper reveals that physical determinism—a basic property that we expect mature physical theories to have or to lack—can require stronger background principles than ZFC if it is understood in the way physics actually uses it. The significance of such results lies in their existence, not just in their frequency or naturalness. This independence holds even at the analytic layer for possible physical systems. And it arises at the regularity layer for the sorts of theoretical questions physicists ask.

*Objection 2. This is only about semantics. The world does not change with axioms.*

The objection is false of the analytic layer and misleading of the regularity layer. First, at the analytic layer, the dependence is counterfactual. Holding fixed a definitional specification of a possible physical system, the coherence, uniqueness, or identity of that system can toggle between V=L and LC/PD. This is as strong as familiar indeterminism thought experiments like Norton's, but of a new kind. It shows that any of three aspects of determinism as studied by philosophers can vary depending on whether one assumes V=L vs. LC/PD in the metatheory.

Second, at the regularity layer, the dependence concerns norms of theory construction, not semantic conventions. Compare general covariance. Its force does not derive from what physicists happen to mean by their words, but from a reasonable physical policy. The world does not care about how we coordinatize it, so neither should our theories of it. Likewise, determinism claims – or giving the opponent the word "determinism" (see Objection 6), "robustness profile" claims – are not bare pointwise assertions. In real physical contexts they are assumed to be typicality laden, and, where nonuniqueness arises, to involve canonicalization. Those commitments force a theory to quantify over rich spaces of admissible representations and to support "almost surely," "generic," or "canonical" verdicts in a representation invariant way. When the resulting determinism profiles lie at $\Sigma^1_2$ complexity, ZFC in general fails to decide whether the regularity conditions required for such verdicts are satisfiable. Once good practice is regimented, ZFC alone can fail to fix what our best physical frameworks say about representation invariant determinism profiles. The fixed presentation Kerr selection problem of Section 6 makes the parallel point at the identity layer.

*Objection 3. Just restrict attention to countably many, or computable, choices.*

The problem with this is that, again, determinism in physics is a fact about the world, not about how we happen to represent it. Whether a system is deterministic should not depend on how limited our descriptions are. Compare again general covariance. If admissible implementations were restricted by pragmatic or epistemic considerations—computability, experimental accessibility, or current modeling practice—then determinism would have to be relativized to those considerations, and would cease to be a representation independent property of the physical system. If determinism is to be an objective feature of the dynamics, one would have to appeal to some physical structure that justifies the restriction to countably many or computable (or whatever) choices, such as a special gauge or a lawlike canonicalization rule. For every physical system where determinism matters, there would have to exist a physically distinguished choice class that is rich enough to capture the intended robustness, tame enough to avoid $\Sigma^1_2$ complexity, and canonically specified by the world. Otherwise we are back to arbitrary choices. But such structure is not supplied by the equations, and introducing it is, as far as I can see, in general unjustified. Restrictions weaken determinism into a practice relative notion of limited philosophical or physical interest. Of course, physicists do not mention all choices, any more than they mention all coordinates. But general covariance demands agreement among all choices of representation, no matter how useless or incomprehensible, and robustness is the same.

*Objection 4. The $\Sigma^1_2$-complete families are engineered. They are not real physics.*

The completeness arguments are explicit embeddings, but they do not depend on unphysical dynamics. Theorem 5.1 uses one fixed nearest neighbor Ising Hamiltonian on $\mathbb{Z}^3$ and one fixed fair zero temperature Glauber update schedule. It is a register construction, not a phase theorem. Theorem 6.2 concerns small continuous Lorentzian perturbations in a fixed collar of a Kerr extension. It does not impose the vacuum Einstein equations, quotient by local isometry, or select among inequivalent continuations. The Kerr result applies to fixed presentation germ canonicalization. The must diverge results show that, once one imposes robustness and typicality norms, $\Sigma^1_2$-hardness is available within mainstream mathematical templates. Of course, many specific applications land below the hardness threshold.  Moving from the present embedded subfamilies to typical disorder laws or to the physically relevant perturbations of a particular black hole is a further reverse-physics question, not something one can safely assume away (Dafermos and Luk 2025; Newman and Stein 1997; Gribov 1978; Singer 1978).

Moreover, for the analytic layer examples, the right standard is not actuality but whether the definitions describe possible systems. The literature on determinism has never required more. Toy models earn their keep by showing what follows from mathematically natural law statements even when no laboratory instantiation is in view. If one now insists that only systems whose specifications are already immune to foundational sensitivity count as “physical,” the conclusion has just been built into the screening criterion.

*Objection 5. You could fix coherence by changing the mathematics (outer measure integrals, non-standard analysis, and so on).*

One can sometimes repair a failure of measurability by adopting alternative conventions, like outer measures, Choquet capacities, or non-standard integration. But that confirms metatheoretic sensitivity, rather than undermining it. Such repairs are interpretive choices, not physical findings. For example, moving from Lebesgue measure to outer measure can restore a form of “measure” for non-measurable sets, but it typically breaks additivity or the link to probability, and it interacts differently with typicality claims. Similarly, adopting a canonical tie break via a constructible well-order yields definable selectors in V=L, but those selectors need not be measurable so do not in general support ensemble reasoning. Determinism practice presupposes an additional regularity theory, and serious axiom candidates disagree about which theory is correct.

*Objection 6. Determinism should be about pointwise well-posedness only. Typicality is extra.*

I do not want to argue about the word “determinism.” If one insists on a very thin notion of determinism (like pointwise existence and uniqueness) then, as I have said, the regularity layer arguments are beside the point of determinism per se. But, whatever you call it, physics often needs more than this. Typicality and robustness are how deterministic structure is used to license predictions, to state stability results, and to articulate claims like SCC, metastability, or robust phase selection. And, to repeat, even if one adopts the thin notion, Section 3 shows that analytic determinism itself is counterfactually metatheory sensitive at all three levels – coherence, uniqueness, and identity – relative to at least “mathematically possible” worlds. While a hard-headed instrumentalist might try to reassure themselves by declaring that nonmeasurable sets are physically impossible, the whole point is that PD/LC and V=L disagree about which sets enjoy that property (and enjoy additional regularity). The insularity thesis fails either way.

*Objection 7. Even if LC/PD and V=L yield different regularity verdicts, nothing observable changes.*

This is incorrect, but irrelevant even if it were correct. At the analytic layer, the dependence on V=L versus LC/PD is directly testable for the possible systems under discussion, in the sense that standard ill-posedness or nonuniqueness results are testable in physics. The regularity layer differs in that the relevant failures manifest as implementation dependence rather than pointwise breakdown. But in experimental and computational practice one can often detect implementation dependence via mesh refinement, cross-code replication, systematic variation of admissible implementations, and the persistence or disappearance of robust verdicts under those changes. Appendix G describes diagnostics of exactly this kind. These can exert methodological pressure on the mathematical background used to underwrite robustness claims. So, whether the mathematics in which a theory is formulated can support the norms that make the theory physically acceptable seems cumulatively testable. Finally, the view that a difference is physically insignificant only if physically observable betrays the crudest kind of positivism whose problems largely defined the philosophy of science in the first half of the 20th century. What physicist with interest in the physical world would claim that the interior of a black hole (Dafermos and Luk 2025) is insignificant?

*Objection 8. We can just accept large cardinals sufficient for LC/PD, and the problem would go away.*

There are two problems with this. First, adopting large cardinals does not eliminate the phenomenon. It relocates it higher in the definability hierarchy. Beyond the Shoenfield absolute level, distinct set theories can disagree not only about regularity properties but even about set membership (Jech 2003; Kanamori 2009; Woodin 1988). Second, this proposal amounts to fixing the standards by which determinism is assessed by mathematical (even set theoretic!) fiat. Questions about whether a physical model is well posed or admits a canonical continuation certainly appear to be matters for physics to settle, not matters to be legislated by logicians.

*Objection 9. Physics uses approximations that are unaffected by these arguments.*

This objection conflates epistemology with metaphysics, like Objection 3. One question is what we can measure or compute. Another is what the theory says the world is like. Determinism is not the claim that predictions are accurate up to some fixed ε, but that the world's evolution is fixed by its present state. Retreating to "deterministic up to ε, for every ε>0" abandons determinism in the ordinary sense unless the ε→0 limit itself is well-defined. But the limit claims physicists rely on – e.g., robust convergence under arbitrarily fine refinement, independence of discretization or windowing in the limit, and so forth – are universal tail statements. Once made precise, and combined with existential quantification over an admissible completion, such claims can reach $\Sigma^1_2$ complexity. Whether such limiting claims are coherent and representation invariant can thus depend on the same set-theoretic regularity principles discussed above.[4]

Are there any other objections that one could give to the arguments of the previous sections? There is one. If such a central concept to science as determinism is sensitive to the V=L vs. LC/PD dichotomy, why has nobody noticed? I suspect that there are three reasons. First, there is the familiar problem of specialization. Descriptive set theorists are not in general physicists, and vice versa, and the projective complexity of physical predicates has not, to my knowledge, been studied in detail. Second, easily foreseeable interactions between set theory and physics (as with non-measurable sets and analytic determinism) seem at most to affect possible physical systems, and so are of limited interest to practicing physicists. Finally, physicists and philosophers have failed to distinguish the two layers of determinism that I have gone to pains to articulate here, analytic and regularity layer determinism. While the analytic layer determinism of actual physical systems may be ZFC robust, regularity layer determinism is not. I suspect that theorists have overlooked this because regularity layer determinism had not been explicitly formulated.

## 9. Reverse physics

Reverse mathematics asks which axioms are needed to prove a mathematical theorem, usually over a weak base theory. Reverse physics asks a related but less uniform question: *which mathematical axioms, regularity principles, and coding assumptions are needed to derive a physical generalization with its intended invariance and typicality content?*

The analogy is deliberately rough. The base theory will often be ZFC rather than a subsystem of second-order arithmetic. The target may be a theorem about a mathematical model rather than a directly observable sentence. And the foundational strength can enter through the semantics of probability or selection, not only through proof-theoretic strength (Friedman 1975; Simpson 2009).

### 9.1 The object of study

Let $T_{phys}$ contain the field equations, constitutive assumptions, boundary or initial conditions, and the explicit idealizations of a physical framework. Let M be a foundational metatheory. Write

$$T_{phys} \vDash_M \Phi$$

when M proves that every intended mathematical model of $T_{phys}$ satisfies Φ. A reverse-physics analysis should record at least five things:

(1) the standard Borel or Polish coding of the states and auxiliary objects;

(2) the exact event, relation, or selector that represents the physical claim;

(3) the probability measure, topology, gauge group, or equivalence relation relative to which regularity is required;

(4) an upper bound on definitional complexity and, where possible, a matching lower bound;

(5) the weakest additional principles sufficient for the desired measurability, uniformization, or absoluteness result.

The fifth item is not automatically a large cardinal question. Compactness, determinacy for a restricted pointclass, countable choice principles, effective uniformization, or a domain-specific uniqueness theorem may suffice.

### 9.2 A worked Kerr agenda

The Kerr discussion separates three targets.

First is *bare extension existence* at a regularity threshold. Under a Borel certificate presentation this is analytic, and Theorem 6.1 rules out exact calibration to the full L[a]-axis in an anti-L[a] universe. The reverse physics question is therefore not "which large cardinal decides bare extendibility?" but whether the intended weak solution or extension relation has the stated Borel presentation and what absoluteness follows from it.

Second is *canonical representation* of an extension germ. Theorem 6.2 proves an $E_0$ obstruction in a fixed collar. Here the metatheoretic split is that V=L supplies a projective selector, while PD excludes one, and ZFC excludes every *universally measurable* selector. The next task is to replace the fixed chart equivalence by local isometry over the common past. That is a rigidity problem, not a routine rewriting.

Third is *selection among inequivalent continuations*. This is the question relevant to recovery of a future at the $C^0$ threshold. The present theorem does not address it. A satisfactory formulation would specify which continuations count as solutions, what equivalence identifies physically identical extensions, and what lawlike regularity a selector must possess. Only then can its projective complexity be classified.

For SCC itself, generic inextendibility and canonical continuation must not be conflated. If a stronger regularity SCC theorem holds, no selector is needed at that threshold because there are no extensions. At $C^0$, where extensions may exist, a selection question arises only if one insists that determinism survive despite extendibility. Reverse physics maps these claims separately.

### 9.3 A statistical mechanical agenda

Theorem 5.1 provides a lower bound inside a fixed local dynamics. The physically important questions concern restrictions.  One is measure. Does a standard Gibbs, metastate, or disorder measure assign full measure to a Borel safe region of the model space? If so, the universal nonmeasurability of the full profile may be irrelevant to the theory's chances. If not, can nonmeasurability be proved for that distinguished measure rather than for an adversarial witness measure?

A second is preparation. Can the existential completion variable be interpreted as microscopic slack in one physically specified macrostate, rather than as an arbitrary infinite code? This requires a mathematically explicit coarse graining map and a physically motivated fiber.

A third is policy. Which boundary conditions, volume sequences, or update schedules are admissible? The full Polish space of coherent paths used in Theorem 5.1 is designed to realize well-foundedness. Standard thermodynamic practice may privilege a smaller class. The restriction must be stated before complexity is assessed.

A fourth is naturality. One should analyze concrete systems—Edwards-Anderson, Sherrington-Kirkpatrick, kinetically constrained models, or metastate constructions—without first planting a universal tree. A negative result showing that their canonical profiles are Borel would also be progress.

### 9.4 Computability and other foundations

A profile may be measurable but lack a *computable* selector. Conversely, restricting to computable policies may lower projective complexity while changing the physical claim. Work on computability in analysis and physics already shows that computable data can generate noncomputable solution behavior in familiar equations (Pour-El and Richards 1989). Reverse physics should keep effective and set-theoretic obstructions separate.

Indeed, the V=L/PD contrast is only one axis. Constructive, predicative, topos-theoretic, and ultrafinitist frameworks can alter which state spaces, limits, or probability assignments are available. The right comparison depends on the physical claim. The present paper uses V=L and PD because they give a clean low projective regularity contrast familiar to mainstream logicians and set theorists.  I do not assume that ZFC occupies a non-negotiable position as foundation.

## 10. Conclusion: determinism after foundations

Physicists and philosophers of physics commonly assume that disagreements in the foundations of mathematics—such as whether there are large cardinals or whether V=L—have no bearing on physical theory (Feferman 1998). This is the insularity thesis, and I have argued that it is suspect.

There are two layers of dependence of determinism on set theory that have not been properly distinguished. At the *analytic* layer, one can construct natural models of possible systems whose determinism profile (coherence, uniqueness, and even the identity of the "unique" solution picked out by a definition) differs between V=L and large cardinal assumptions sufficient to imply PD. This already shows that the *concept* of determinism is entangled with higher set theory. At the *regularity* layer, determinism demands invariance across admissible implementations, and is expressed in typicality and canonicalization language. Making such claims precise and representation invariant requires regularity and selection conditions for definable sets and relations. Those conditions are generally independent of ZFC and divide between V=L and LC/PD. The results established here include that a fixed Ising Hamiltonian yields a $\Sigma^1_2$-complete tail-readout coding theorem, and a fixed Kerr collar yields an $E_0$ obstruction to regular selection of representative metric-germ codes. The first is not a theorem about a standard phase or disorder measure. The second does is not an SCC theorem. Deeper physical questions belong to reverse physics (Earman 1986; Loewer 2020; Dafermos and Luk 2025).

Let me emphasize that the point has not been the obvious one that there are seemingly physical questions (e.g., about the behavior of computers searching for a contradiction in PA) that one can only decide by strengthening the axioms of mathematics (by, say, adding the assumption of their consistency). That is just Gödel's theorem (Gödel 1931), assuming that the relevant principles are consistent and recursively axiomatizable (and represent all primitive recursive functions). V=L and LC are seriously entertained extensions of standard mathematics (Maddy 1997; Steel 2014; Koellner 2009), and the undecidables in question are in no way contrived. We have to take a stand on V=L and LC/PD like the stand that the pioneers of foundational studies took with respect to the Axiom of Choice (Jech 2003, Ch. 1; Kanamori 2009, §0). The question of which of V=L or LC/PD is true -- or whether it even makes sense to say that one of them is true -- is not a question that admits of proof. It is a question, to borrow a turn of phrase, of "deeply rooted mathematical taste". As Jensen writes:

> Most proponents of V=L and similar axioms support their belief with a mild version of Ockham's razor. L is adequate for all of mathematics. It gives clear answers to deep questions. It leads to interesting mathematics. Why should one assume more? … I do not understand … why a belief in the objective existence of sets obligates one to seek ever stronger existence postulates. Why isn't Platonism compatible with the mild form of Ockham's razor? … I doubt that one could, with the sort of evidence I have, convert the mathematical world to one or the other point of view. Deeply rooted differences in mathematical taste are too strong and would persist. (Jensen 1995, 400–401)

One reaction to this situation is that, despite widespread opinion to the contrary, Quine was right, and the search for new axioms is not different from the search for new physical laws (Quine 1951, 1990). As he puts it, "sentences such as the continuum hypothesis…which are independent of [standard] axioms, can…be submitted to the considerations of simplicity, economy, and naturalness that contribute to the molding of scientific theories generally. Such considerations support… [the] Axiom of Constructibility, 'V=L'" (Quine 1990, 95). This position raises familiar quandaries. Set theorists are not responsive to experiment in the way physicists are (Maddy 1997, 155; 2011). But this much is true: when we say that general relativity predicts X, we mean that X follows from the field equations plus boundary conditions, where "follows from" is a relationship of entailment within a mathematical framework like ZFC + LC/PD. If that entailment depends on which set theoretic axioms we assume, then the theory's content is metatheory sensitive. This is not a defect of our measurements or approximations. It is a fact of theorizing. Recognizing this clarifies the relationship between mathematical and physical reasoning, showing that they are allied in a way that non-Quineans have not appreciated.

My own view (Clarke-Doane 2022, 2024, Forthcoming) is that no sense can be made of the claim that either V=L or LC/PD is "really right." Each affords an arena in which to carry out mathematical reasoning, broadly like Euclidean and non-Euclidean geometries. The heady difference is that all the geometries can live in a single set theory. But any set theory takes itself to be maximal, to govern everything. This raises unresolved questions about how to formulate the "monism–pluralism" debate, and whether it is even factual. If one is a set-theoretic pluralist, however, then, by the arguments given here, one must be a pluralist about core physical concepts too. ("Concepts" rather than "concept" because if determinism is relative to set-theoretic foundations, then so too, presumably, is physical necessity, prediction, explanation, causation, and more.) As Quine writes in a different context, "Carnap has maintained that [the question of which axioms hold] is a question not of matters of fact but of choosing a convenient language form, a convenient conceptual scheme or framework for science. With this I agree, but only on the proviso that the same be conceded regarding scientific hypotheses generally" (1951, 43).[5]

## Appendix A. Analytic-layer constructions

### A.1 The definable Fubini pathology in L

Assume V=L. There is a parameter-free lightface $\Delta^1_2$ well-order $\prec$ of $2^\omega$ whose initial segments are countable. Let

$$W=\{(x,y)\in 2^\omega\times 2^\omega : x\prec y\}.$$

Then W is lightface $\Delta^1_2$. For every fixed y, the horizontal section

$$W^y=\{x : x\prec y\}$$

is countable. For every fixed x, the vertical section

$$W_x=\{y : x\prec y\}$$

is co-countable. If W were measurable for the product of the fair-coin measure with itself, Fubini's theorem would give product measure 0 from the horizontal sections and product measure 1 from the vertical sections. Hence W is not Lebesgue measurable and therefore not universally measurable. The construction is a definable refinement of Sierpiński's classical examples (Sierpiński 1924).

It also lacks the Baire property.. The Kuratowski–Ulam theorem states that a set with the Baire property in a product of second countable Baire spaces is meagre if and only if comeagrely many of its vertical sections are meagre (Kechris 1995, Section 8.K); the coordinate-swap homeomorphism transfers the same criterion to horizontal sections. Suppose W had the Baire property. Every horizontal section $W^y = \{x : x \prec y\}$ is a countable initial segment, hence meagre. By Kuratowski–Ulam, W is meagre. Now use the other family. Every vertical section $W_x = \{y : x \prec y\}$ is co-countable, so every vertical section of the complement $W^c$ is countable, hence meagre. Since $W^c$ has the Baire property as well, Kuratowski–Ulam makes $W^c$ meagre; equivalently, W is comeagre. No set is both meagre and comeagre in the Baire space $2^\omega \times 2^\omega$. Hence W does not have the Baire property in L.

Transporting W through fixed Borel isomorphisms gives the version on $(0,1)^2$ used in the main text. Under PD every projective set has all three regularity properties, so the same defining scheme behaves differently in the two metatheories.

**A.2 Coherence of the weak heat equation**

Let $\Omega=(0,1)^2$ and transport the preceding W to a subset of Ω. Put $V=1_W$ and consider

$$\partial\square u=\Delta u-Vu,\quad u(0,\cdot)=u_0,$$

with the intended weak form on $H^1_0(\Omega)$

$$a(u,v)=\int_\Omega(\nabla u\cdot\nabla v+Vuv)\,dx.$$

Let

$$b(x_1,x_2)=x_1(1-x_1)x_2(1-x_2).$$

Then $b\in H^1_0(\Omega)$ and $b>0$ in Ω. If $b^2 1_W$ were measurable, then

$$W=\{x\in\Omega : b^2(x)1_W(x)>0\}$$

would be measurable. Thus in L the integrand in a(b,b) is not Lebesgue measurable, and the intended bilinear form is not defined. Under LC/PD, V is bounded and measurable, so the standard weak formulation and semigroup theory apply (Brezis 2011; Evans 2010).

The conclusion concerns this formulation. It does not preclude replacing the integral, changing the coefficient class, or declaring the defining formula physically inadmissible.

### A.3 Uniqueness by a convexity switch

Let $U\subseteq\mathbb{N}\times\mathbb{R}$ be universal for the lightface $\Pi^1_2$ subsets of $\mathbb{R}$, and put

$$\psi\Leftrightarrow\forall\ n\ (U_n \text{ is Lebesgue measurable}).$$

Define $\gamma=-1$ if $\psi$ and $\gamma=+1$ otherwise. Under LC/PD, $\psi$ holds. In L, the lightface $\Delta^1_2$ nonmeasurable set of A.1 is a $\Pi^1_2$ section, so $\psi$ fails.

For

$$E_\gamma(z)=z^4-\gamma\ z^2,$$

if $\gamma=-1$, then $E''_\gamma(z)=12z^2+2>0$ and the unique minimizer is 0. If $\gamma=+1$, the global minimizers are $\pm1/\sqrt{2}$. On a Hilbert space H with unit vector e, the functional

$$\mathscr{E}_\gamma(u)=E_\gamma(\langle u,e\rangle)+\|u-\langle u,e\rangle e\|^2$$

has the corresponding one or two minimizers. The variational calculation is elementary; all metatheoretic work is in the definition of $\gamma$.

### A.4 Smooth identity shift

Let

$$C=\{n:U_n \text{ is not Lebesgue measurable}\}.$$

Choose nonzero $C^\infty$ bumps $\varphi_n$ with pairwise disjoint supports and

$$\|\varphi_n\|_{C^m}\leq2^{-n}\ \text{ whenever } m\leq n.$$

For each fixed m, the tail from $n\geq m$ is summable in the $C^m$ norm. Hence

$$u_0=\Sigma_{n=0}^{\infty}2^{-n}1_C(n)\varphi_n$$

converges in every $C^m$ norm and therefore in $C^\infty$. Under LC/PD it vanishes. In L it is nonzero because C is nonempty and the supports are disjoint. Standard continuous dependence for a well-posed linear equation then gives different unique evolutions from the same intensional specification (Hadamard 1923; Evans 2010).

## Appendix B. Regularity and preservation facts

### B.1 Universal measurability

For a Polish space X, a set $A\subseteq X$ is universally measurable when it belongs to the completion of every Borel probability measure on X. The definition is tied to the standard Borel structure, not to one topology among all topologies inducing it.

If $f:Z\rightarrow X$ is Borel and A is universally measurable, then $f^{-1}(A)$ is universally measurable. Given a Borel probability $\mu$ on Z, let $\nu=f_*\mu$. There are Borel sets $B_0\subseteq A\subseteq B_1$ with $\nu(B_1\setminus B_0)=0$. Then

$$f^{-1}(B_0)\subseteq f^{-1}(A)\subseteq f^{-1}(B_1)$$

and $\mu(f^{-1}(B_1\setminus B_0))=0$. Thus the preimage belongs to the completion of $\mu$.

This preservation fact is what allows continuous completeness reductions to transmit universal-measurability failure.

### B.2 The Baire property

A set $A \subseteq X$ has the Baire property when $A \triangle O$ is meagre for some open O. Unlike universal measurability, the Baire property is not preserved under arbitrary continuous preimages. If a continuous map has a large constant fiber lying inside a meagre exceptional set, its preimage can be arbitrary.

Concretely, let $A \subseteq 2^\omega$ be $\Sigma^1_2$-complete. Embed $2^\omega$ homeomorphically as a closed nowhere dense subset C of another copy of $2^\omega$ and let A' be the image of A. Then A' is still $\Sigma^1_2$-complete: compose reductions with the embedding. But $A' \subseteq C$ is meagre and therefore has the Baire property. Completeness alone cannot yield Baire failure.

### B.3 Projective regularity under PD

Under PD, every projective set of reals is Lebesgue measurable, universally measurable, and has the Baire property. These facts transfer to projective subsets of arbitrary Polish spaces through standard Borel presentations. PD also yields strong uniformization and scale results, but the paper invokes only the stated regularity consequences except in the comparison with projective selector machinery (Kechris 1995; Moschovakis 2009; Martin and Steel 1989).

### B.4 Uniformization, selector maps, and transversals

For a relation $R \subseteq X \times Y$, a uniformizing function f satisfies $(x,f(x)) \in R$ on the projection of R. For an equivalence relation E on X, a selector map in the sense used in Section 6 is constant on each class and chooses an E-equivalent representative. Its range is a transversal, a set meeting every class exactly once.

These notions should be distinguished from a selector that merely picks some y in each vertical section of a general relation. The no-go theorem for $E_0$ is about class constant selector maps and transversals. Ordinary choice gives such a map, but no universally measurable or Baire measurable one exists.

## Appendix C. The universal-tail classification

Let $X,Y,\mathbb{P}$ be Polish and let

$$G \subseteq X \times Y \times \mathbb{P} \times \mathbb{N} \times \mathbb{N}$$

be Borel. Write $G(x,y,\rho,k,n)$ for membership. Define

$$B(x,y,\rho) \Leftrightarrow \forall k \ \exists N \ \forall n \geq N \ G(x,y,\rho,k,n).$$

Because all three quantifiers here range over $\mathbb{N}$, B is Borel. Its precise Borel rank may be high and is irrelevant. Now

$$UT(x,y) \Leftrightarrow \forall \rho \in \mathbb{P} \ B(x,y,\rho).$$

The complement is the projection of the Borel set

$$\{(x,y,\rho): \neg B(x,y,\rho)\},$$

so UT is coanalytic. Finally,

$$\mathrm{DET}(x) \Leftrightarrow \exists\, y \in Y\ \mathrm{UT}(x,y)$$

is $\Sigma^1_2$.

Suppose DET is complete for that pointclass under continuous reductions. Internally to V=L, let D be the lightface $\Delta^1_2$ non-universally-measurable set from Appendix A, transported to $2^\omega$. Completeness gives a continuous f with

$$x \in D \Leftrightarrow f(x) \in \mathrm{DET}.$$

If DET were universally measurable, Appendix B.1 would make D universally measurable. Hence it is not. Under LC/PD, projective regularity gives universal measurability and the Baire property directly.

## Appendix D. The fixed-Hamiltonian Ising register

### D.1 Tree normal form

Let Tr be the Polish space of trees on $\omega^{<\omega}$, viewed as a closed subspace of $2^{(\omega^{<\omega})}$. The set

$$\mathrm{WF}=\{T \in \mathrm{Tr}: T \text{ is well founded}\}$$

is $\Pi^1_1$-complete. The usual tree normal form for $\Sigma^1_2$ sets supplies a continuous map

$$(\alpha,\beta) \mapsto T_{\alpha,\beta}$$

such that

$$\mathrm{WF}_{\exists} =\{\alpha: \exists\, \beta\ T_{\alpha,\beta} \in \mathrm{WF}\}$$

is $\Sigma^1_2$-complete. Continuity means that for every finite node s, the truth of $s \in T_{\alpha,\beta}$ is determined by finite initial segments of α and β.

### D.2 One fixed Hamiltonian

Partition $\mathbb{Z}^3$ into dimers

$$D_{j,k,\ell}=\{(2j,k,\ell),(2j+1,k,\ell)\}.$$

Set $J_{uv}=1$ when $\{u,v\}=D_{j,k,\ell}$ for some triple and $J_{uv}=0$ for every other nearest-neighbor pair. Let the external field vanish. The local Hamiltonian difference at a site v depends only on its paired neighbor $v^*$:

$$\Delta_v H=2\sigma(v)\sigma(v^*).$$

The unique zero-temperature minimizing update is therefore

$$\sigma(v) \leftarrow \sigma(v^*).$$

Choose a computable schedule $q:\mathbb{N}\to\mathbb{Z}^3$ such that each site occurs infinitely often. This fixes the evolution law completely. Every configuration constant on each dimer is stationary under q.

Fix a computable injection $c:\omega^{<\omega}\to\mathbb{Z}^3$ whose range indexes dimers rather than individual sites; write $D_s$ for the dimer assigned to s. Given (α,β), define $\sigma_{\alpha,\beta}$ by

$$\sigma_{\alpha,\beta}\restriction D_s = (-1,-1) \text{ if } s \in T_{\alpha,\beta}, \text{ and } \sigma_{\alpha,\beta}\restriction D_s = (+1,+1) \text{ if } s \notin T_{\alpha,\beta}.$$

and put all remaining dimers in the (+1,+1) state. Each output coordinate depends on one clopen tree-membership fact, so $(\alpha,\beta)\mapsto\sigma_{\alpha,\beta}$ is continuous. Every such state is stationary.

### D.3 The path space and verdict

Let

$$\mathbb{P}= \{\rho\in(\omega^{<\omega})^{\omega}: \rho(n)\subsetneq\rho(n+1) \text{ for every } n\}.$$

Failure is witnessed at a finite stage, so $\mathbb{P}$ is closed. For $\rho\in\mathbb{P}$ define $V(\alpha,\beta,\rho,n)$ to be the common spin on $D_{\rho(n)}$ after stage n of the fixed evolution. Stationarity gives

$$V(\alpha,\beta,\rho,n)=+1 \Leftrightarrow \rho(n)\notin T_{\alpha,\beta}.$$

The relation is Borel, indeed continuous when the discrete verdict space is used.

### D.4 Exact decoding

If $T_{\alpha,\beta}$ has an infinite branch b, let $\rho(n)=b\restriction(n+1)$. Then $V=-1$ at every stage. Conversely, if $T_{\alpha,\beta}$ is well founded, every strict chain leaves it at some stage. Once a chain leaves a tree, every later extension remains outside. Hence

$$T_{\alpha,\beta}\in WF \Leftrightarrow \forall\rho\in\mathbb{P}\ \exists\ N\ \forall\ n\geq N\ [V(\alpha,\beta,\rho,n)=+1].$$

Projecting over β proves Theorem 5.1. Notice that no parameter enters the Hamiltonian, the coupling pattern, the update schedule, or the verdict rule.

### D.5 Why the construction is not a phase theorem

The spin read at stage n is the spin on a newly designated dimer, not the magnetization of an expanding physical region. The path ρ is a tree following diagnostic, not an ordinary asynchronous update order. The exact microstate is stationary, so no phase transition or relaxation phenomenon is being encoded.

## Appendix E. Bare extension and the Shoenfield obstruction

### E.1 Coding assumption

Let $\mathscr{D}$ be a standard Borel data space and $\mathscr{Y}$ a Polish certificate space. Fix a real a coding the presentation. Assume

$$R(d,y,a)$$

is Borel and expresses the local compatibility, regularity, boundary matching, and equation requirements imposed on a candidate extension. Then

$$\mathrm{Ext}(d) \Leftrightarrow \exists\ y\ R(d,y,a)$$

is $\Sigma^1_1(a)$.

This appendix does not prove that every geometric extension predicate admits such a presentation. For classical $C^k$ extensions on fixed compact charts, Borel coding is standard. For weak solution concepts, quotient presentations, or maximal developments constructed only indirectly, that point requires separate analysis.

### E.2 Exact calibration collapses the real universe

Let $F:\omega^{\omega}\to\mathscr{D}$ be Borel with code absorbed into a. Then

$$A=\{x:\mathrm{Ext}(F(x))\}$$

is $\Sigma^1_1(a)$. Choose an arithmetical formula $\varphi$ such that

$$x\in A\Leftrightarrow \exists z\, \varphi(x,z,a).$$

Assume $A=\mathbb{R}\cap L[a]$.

For every $x\in L[a]$, the universe satisfies $\exists z\, \varphi(x,z,a)$. By Shoenfield absoluteness, so does L[a]. Therefore

$$L[a]\vDash\forall x\, \exists z\, \varphi(x,z,a).$$

This is a $\Pi^1_2(a)$ sentence. A second application of Shoenfield absoluteness yields

$$V\vDash\forall x\, \exists z\, \varphi(x,z,a).$$

Thus $A=\mathbb{R}$, and the assumed equality gives $\mathbb{R}\subseteq L[a]$.

The theorem does not prove V=L[a]. Sets above the reals may remain outside L[a]. The contradiction needed by the application is simply the existence of a non-L[a] real.

### E.3 Consequences for SCC coding programs

Suppose a proposed Kerr construction has three ingredients:

(1) a Borel map from real codes to initial data;

(2) a bare extension predicate presented analytically;

(3) exact equivalence of bare extension with membership in L[a].

In any anti-L[a] universe the three cannot coexist. Strengthening a local screen until it becomes equivalent to bare extension does not evade the theorem but transfers the analytic bound to the screen predicate. A positive construction must change at least one ingredient. It may target a higher complexity adversarial or canonicalization predicate, use a non-Borel data map, weaken exact calibration, or work in a universe with no non-L[a] reals.

## Appendix F. The fixed-collar $E_0$ theorem

### F.1 The perturbation space

Let $N=K\times[0,1]$, with K a compact manifold with boundary if necessary, and fix a continuous Lorentzian metric $\bar{g}_0$ on N. Let $H=K\times\{0\}$ be the horizon boundary. Choose a continuous background Riemannian metric to define tensor norms and put

$$B=C_0(N;\mathrm{Sym}^2T^*N) =\{h\in C(N;\mathrm{Sym}^2T^*N):h|_H=0\}.$$

This is a separable Banach space in the uniform norm. Since $\bar{g}_0$ is nondegenerate and N is compact, there is $\varepsilon>0$ such that $\bar{g}_0+h$ has Lorentzian signature for every h in the closed ball

$$X=\{h\in B:\|h\|_\infty\leq\varepsilon\}.$$

Regard X as the code space. Adding $\bar{g}_0$ converts a code into a continuous Lorentzian extension with the same boundary values on H.

For $m\geq1$ let

$$C_m=\{(h,h')\in X^2:h=h' \text{ on } K\times[0,1/m]\}.$$

Each $C_m$ is closed. Hence

$$E_{germ}=\cup_{m\geq 1}C_m$$

is $F_\sigma$. It is plainly an equivalence relation.

**F.2 Marker embedding**

Choose points $(p_n,s_n)$ with $s_n\downarrow 0$ and pairwise disjoint coordinate balls $B_n$ whose closures lie in $K\times(0,1]$. Choose the balls so that the upper s-coordinate of $B_n$ tends to 0. Choose continuous, or smooth in the interior, symmetric tensor bumps $u_n$ supported in $B_n$, with $u_n(p_n,s_n)\neq 0$. Rescale them so that

$$\|u_n\|_\infty\leq 2^{-n-2}\varepsilon.$$

For $\gamma\in 2^\omega$ put

$$\Theta(\gamma)=\Sigma_{n=0}^{\infty}\gamma(n)u_n.$$

The supports are disjoint, and the norms tend to zero, so the sum converges uniformly and vanishes continuously at H. It lies in X. If two sequences agree through N-1, then

$$\|\Theta(\gamma)-\Theta(\gamma')\|_\infty\leq \sup_{n\geq N}\|u_n\|_\infty,$$

so $\Theta$ is continuous. It is injective because the value in $B_n$ reads the nth bit. Since $2^\omega$ is compact and X is Hausdorff, $\Theta$ is a homeomorphism onto its compact image F.

If $\gamma\ E_0\ \gamma'$, the two sums differ in only finitely many marker balls. A sufficiently small collar excludes all of them, so the perturbations define the same germ. Conversely, if the sequences differ at infinitely many indices, there are differing marker balls in every collar, so the germs are unequal. Therefore

$$\Theta(\gamma)\ E_{germ}\ \Theta(\gamma')\Leftrightarrow\gamma\ E_0\ \gamma'.$$

**F.3 No universally measurable selector**

Let Y be either X or F, with its subspace Borel structure. Suppose

$$S:Y\rightarrow X$$

is a selector map: $S(h)\ E_{germ}\ h$ for every $h\in Y$, and $S(h)=S(h')$ whenever $h\ E_{germ}\ h'$. Assume S is universally measurable. Then

$$r=S\circ\Theta:2^\omega\rightarrow X$$

is measurable for the fair Bernoulli product measure $\mu$ and is invariant under $E_0$.

Let $\Gamma$ be the countable group of finite coordinate flips. Its orbits are the $E_0$ classes. If $A\subseteq 2^\omega$ is $\mu$-measurable and $\Gamma$-invariant, then $\mu(A)\in\{0,1\}$. Indeed, for every n, invariance makes $\mu(A\cap[w])$ the same for all cylinders [w] of length n. Summing over those cylinders gives

$$\mu(A\cap[w])=2^{-n}\mu(A)=\mu([w])\mu(A).$$

Thus A is independent of the algebra of cylinders and hence of the Borel $\sigma$-algebra it generates. Since A is measurable only for the completion, choose a Borel set B with $\mu(A\triangle B)=0$. Independence applies to B, and $\mu(A\cap B)=\mu(A)$ while $\mu(B)=\mu(A)$, so $\mu(A)=\mu(A\cap B)=\mu(A)\mu(B)=\mu(A)^2$.

For every Borel $C\subseteq X$, the set $r^{-1}(C)$ is invariant and measurable, so the pushforward $\nu=r_*\mu$ takes only the values 0 and 1 on Borel sets. A 0-1 Borel probability on a Polish space is a point mass.

One proof uses tightness. Choose a compact $K_0$ with $\nu(K_0)=1$. Cover $K_0$ by finitely many sets of diameter below 1 and form a finite Borel partition subordinate to the cover. One piece has measure 1; its closure inside $K_0$ is compact and has small diameter. Repeat inside that compact set with diameters tending to zero. Continuity from above gives measure 1 to the unique point in the nested intersection.

Hence $r(\gamma)=p$ for $\mu$-almost every $\gamma$. For every such $\gamma$,

$$p=S(\Theta(\gamma))\ E_{germ}\ \Theta(\gamma).$$

Any two points in this conull set therefore have $E_0$-equivalent codes, by transitivity and the marker equivalence. But one $E_0$ class is countable and $\mu$-null. Contradiction.

**F.4 No Baire measurable selector on F**

Suppose $S:F\to X$ is Baire measurable and a selector. Again let $r=S\circ\Theta$. If $A\subseteq 2^\omega$ is invariant under finite flips and has the Baire property, then A is meagre or comeagre. For if A is nonmeagre, it is comeagre in some basic cylinder [w]. Finite flips carry [w] to every cylinder of the same length and preserve A, so A is comeagre in their union, which is all of $2^\omega$.

Fix a countable base $(U_i)_{i\in\mathbb{N}}$ for X. Each $r^{-1}(U_i)$ is invariant and has the Baire property, so it is meagre or comeagre. Let G be the intersection of all comeagre preimages and all complements of meagre preimages. Then G is comeagre. Any two points of G have images belonging to exactly the same basic open sets, and therefore have equal images. Thus r is constant on G. As before, G would lie in one countable $E_0$ class, impossible.

It remains to rule out a Baire measurable selector on all of X. For $m\geq 1$ let $G_m$ be the closed subspace of perturbations in B vanishing on $K\times[0,1/m]$, and let $G=\cup_m G_m$. Every element of B vanishes at s=0, so G is dense in B. For $k\in G_m$, translation by k is a homeomorphism of B, and it preserves $E_{germ}$ wherever both points lie in X, because h and h+k agree on $K\times[0,1/m]$. Suppose $S:X\to X$ is a Baire measurable selector and fix a countable base $(U_i)$ for X. Each set $A_i=S^{-1}(U_i)$ is $E_{germ}$-invariant and has the Baire property. If $A_i$ is nonmeagre, it is comeagre in some open ball. Given any interior point $h_1$ of X, density of G supplies a translation carrying a small ball inside that ball to a neighborhood of $h_1$ within X, and invariance transports comeagreness there. Since the interior is dense in X, each $A_i$ is meagre or comeagre. Intersecting as before yields a comeagre set on which S is constant, so one $E_{germ}$ class is comeagre. But each class is a countable union of sets $X\cap(h+G_m)$, and each $h+G_m$ is a proper closed affine subspace of B, hence nowhere dense. Every class is therefore meagre. This contradiction rules out the selector on X.

**F.5 The constructible selector**

Assume V=L and fix a real code for the Polish presentation of X. Let $\leq_L$ be the canonical projective well-order of the reals. Define

$$S(h)=\text{the } \leq_L\text{-least } h'\in X \text{ such that } h'\ E_{germ}\ h.$$

This is total and constant on classes. Its graph is defined by

$$\Gamma_S(h,h') \Leftrightarrow h'\ E_{germ}\ h\ \wedge\ \forall z\in X(z\ E_{germ}\ h\to h'\leq_L z).$$

The germ relation is Borel and $\leq_L$ is projective. A further real quantifier leaves the graph projective. The proof makes no claim that the graph lies at $\Delta^1_2$.

By F.3 and F.4, the selector is neither universally measurable nor Baire measurable on the encoded $E_0$ family.

### F.6 PD excludes projective selectors

Assume PD and suppose a selector S on X, or on F, has projective graph. For every open $U \subseteq X$,

$$S^{-1}(U) = \{h : \exists\ y\ (\Gamma_S(h,y) \wedge y \in U)\}$$

is projective. Under PD it is universally measurable and has the Baire property. Hence S is universally measurable and Baire measurable, contrary to F.3 or F.4. This proves the metatheoretic dichotomy.

### F.7 Scope

The proof uses only a continuous Lorentzian collar and equality in a fixed presentation. It applies to exact Kerr, to a Dafermos-Luk collar when the hypotheses of that theorem are met, and to many non-Kerr extension problems. That breadth is mathematically useful and physically limiting. Nothing in the proof detects the Einstein equations beyond the horizon, local isometry, curvature invariants, or inequivalence of futures.

## Appendix G. Empirical and simulation considerations

No finite experiment can survey all reals in a completion or policy space. Nor would a failure of numerical convergence discriminate directly between V=L and PD. The possible empirical relevance is indirect. A physical theory may claim that a verdict is stable under mesh refinement, gauge choice, boundary slack, or reconstruction method. Simulations can test finite families of such choices and increasing resolution. Persistent disagreement can reveal that the stated robustness idealization lacks a canonical implementation. The remedy may be a new physical restriction, a stronger analytic theorem, a different topology, or an explicit convention. Higher set theory is only one possible diagnosis, and usually not the first. Conversely, successful cross-code agreement does not prove universal measurability. It provides evidence that the tested region and tested policies are tame. A reverse physics theorem would explain why—perhaps by proving that the physically distinguished measure is concentrated on a Borel subset, that compactness collapses a real quantifier, or that a canonical gauge reduces the quotient.

The evidential relation is methodological. Mathematics supplies conditions under which ensemble and selection claims are well defined. Physics supplies reasons to accept one state space, measure, topology, or equivalence relation rather than another. Neither side can be replaced by the slogan that axioms are empirically tested, or by the opposite slogan that foundations are irrelevant.

## Endnotes

[1] Thanks to Avner Ash for catching a mistake in my treatment of PDEs; Will Cavendish for pressing me to consider the discrete case; Gabriel Goldberg for illuminating discussion of complexity classes and projective uniformization in L; Joel David Hamkins for general exchanges about the relevance of incompleteness to determinism; Douglas Blue for the virtues of LC/PD; John Norton and William D’Alessandro for feedback on a talk I gave on the topic at UNC Chapel Hill; Peter Galison, Daine L. Danielson, Siddharth Muthukrishnan, and attendees at my Harvard Foundations Black Hole Initiative talk on the subject; two anonymous referees for helpful reports; and the students in my Undecidability in Physics Fall 2025 graduate seminar at Columbia for probing questions about the project. This article is an experiment for me. I have been writing it over the course of or a couple of years in a vigorous back-and-forth with AI systems. I used Claude (Anthropic) and GPT (OpenAI) to help me think across fields, test arguments, track down references, construct proofs, and sometimes to suggest presentation. They

have made it possible for me to work across descriptive set theory, mathematical relativity, engineering, quantum foundations, and related areas in a way I could not have managed on my own. The claims of the paper, its structure, and, of course, any mistakes are my responsibility. (To my knowledge, nothing has been written on the interaction between V=L/PD and determinism. But for important precursors to some of the general issues discussed here, see Farah & Magidor 2012 and Pour-El & Richards 1989 and the references therein.)

[2] For discussion of other “restrictive” axioms, besides V=L, see Fraenkel, Bar-Hillel & Lévy (1973, §6.4). The dominant narrative among set theorists who are willing to speak of truth or falsity is that V=L is false, and “clearly” so. Maddy (1997, Pt II, § 4) contains a nice explication of the reasons advocates of this narrative supply. But this view is “not unanimously shared” (Fontanella 2019, 32). Devlin writes, “What is my own view?...Currently I tend to favour [V=L]....At the moment I think I am in the majority of informed mathematicians, but the minority of set theorists...” (1981, 205). Arrigoni says, “I believe it perfectly in order to characterize . . . ZFC + V=L as intuitively plausible . . . ” (2011, 355). Jensen writes, “I personally find it a very attractive axiom” (1995, 398), and confesses “to being emotionally a Pythagorean” (1995, 401, n. 17), Pythagoreanism being his name for the orientation on which V=L embodies “an entirely coherent Pythagorean picture of the world” (1995, 401). And Pinter writes, “there is a strong intuitive basis for considering L to be the class of all sets. By definition, L contains all the sets that are describable by a formula in the language of set theory. And there is no practical reason to admit sets which lack any description, for we would never make use of such sets. They would merely sit there and muddy the waters. Thus, from this point onward we shall assume the following important axiom: Axiom of Constructibility: Every set is constructible, that is, every set is in L. This axiom is usually denoted by the symbol V=L” (2014, 227). We will meet another supporter in the Conclusion.

[3] Note that for each fixed $(\alpha,\beta,\varrho)$, the evolution is deterministic. The determinism claim here is the robust verdict “the refinement path phase diagnosis eventually locks to +,” which is (i) a universal-tail robustness requirement over a rich policy space $\mathbb{P}$ and (ii) quantified over completions $\beta$. This is the recipe that forces $\Sigma^1_2$ complexity and hence activates the must diverge route.

[4] For each fixed tolerance $\varepsilon>0$, a claim of the form “the diagnostic succeeds to accuracy $\varepsilon$ beyond some finite refinement stage” is tame. But the robustness claims used in physics concern the $\varepsilon\to 0$ limit. Formally, let $\mathbb{P}$ be a Polish space of admissible policies and let Good$(x,\varrho,k,n)$ be a Borel predicate meaning that, at refinement stage n under policy $\varrho$, the diagnostic succeeds to tolerance $\varepsilon_k$, where $\varepsilon_k\downarrow 0$. The $\varepsilon\to 0$ robustness claim has the form: $\forall k\in\mathbb{N}\ \forall \varrho\in\mathbb{P}\ \exists N\ \forall n\geq N$ Good$(x,\varrho,k,n)$. This is a universal tail condition. For each policy and each tolerance, success eventually stabilizes. With Good Borel and $\mathbb{P}$ Polish, the resulting predicate is $\Pi^1_1$ in x. If one then forms a determinism profile by existentially quantifying over admissible completions y (as in the main text), one obtains a set of the form: DET$(x) \Leftrightarrow\ \exists y\ \forall k\ \forall \varrho\ \exists N\ \forall n\geq N$ Good$(x,y,\varrho,k,n)$, which is $\Sigma^1_2$. So, even when each fixed-$\varepsilon$ approximation claim is ZFC-absolute, the $\varepsilon\to 0$ robustness claim lives at $\Sigma^1_2$ complexity. By Theorem 4.3, such limit claims can be universally measurable under LC/PD while failing universal measurability under V=L. Finite precision practice does not insulate determinism from metatheory. It actually presupposes limiting claims whose status ZFC alone need not settle.

[5] Here is the idea: any careful statement of either view must be made against some background theory that fixes what counts as a universe of sets. Inside that framework, “V” denotes every set

that the framework recognizes, so the monist's claim that there is one true V is trivially true by that framework's own lights. But taken independently of any background the claim has no content. One can choose to interpret it using another theory. But different choices yield different contents. The same holds for pluralism. Its "many universes" are sets or proper classes in one background, not an external plurality that both sides could recognize. What remains is, if not a verbal question about the meaning of the word "set" among select academics, a pragmatic question—which framework to use—not a factual one about what exists (Clarke-Doane 2022, Ch. 4; Forthcoming). Of course I do not pretend to have properly defended that view here.

[6] The philosophical pattern emphasized in the main text is not altogether peculiar to the paper's $\Sigma^1_2$-completeness constructions. Levy (2025) proves that universally or Lebesgue measurable equilibrium selection for Borel parametrized games is independent of ZFC and equivalent to the universal measurability of all $\Delta^1_2$ sets, even in the two-player zero sum case. Backs, Jackson, and Mauldin (2015) show that Maharam's question about uniformly σ-finite disintegrations can depend on set theory, and recast nonuniformity in terms of whether the atom relation is a countable union of Borel graphs. Carassus and Ferhoune (2024, 2025) show that dynamic optimization under model uncertainty uses PD to recover measurable/projective selector machinery and ε-optimal selectors, while also giving an L-side example in which the unique optimal strategy is not Lebesgue measurable. Again, these are not cases of physical determinism. But they show that the contrast between pointwise existence and regular global selection already appears, without bespoke completeness coding, in natural theories people actually study.